\documentclass[11pt]{article}
\usepackage{amsfonts}
\usepackage{amssymb}
\usepackage{amsfonts,epsfig}
\begin{document}

\newtheorem{them}{Theorem}
\newtheorem{lema}{Lemma}
\newtheorem{prop}{Proposition}
\newtheorem{definition}{Definition}
\newtheorem{rema}{Remark}
\newtheorem{coro}{Corollary}
\newtheorem{exam}{Example}

\newcommand{\ssq}{$\square$}
\newcommand{\cl}{{\cal L}}
\newcommand{\qed}{\hfill$\Box$}
\newcommand{\pr}{{\bf Proof. }}
\newcommand{\rr}{\Bbb{R}}
\newcommand{\cc}{\Bbb{C}}
\newcommand{\wh}{\widehat}
\newcommand{\ra}{\Rightarrow}
\newcommand{\al}{\alpha}
\newcommand{\be}{\beta}
\newcommand{\la}{\lambda}
\newcommand{\lab}{\lambda_{\alpha, \beta}}
\newcommand{\uab}{\upsilon_{\alpha, \beta}}
\newcommand{\oab}{\omega_{\alpha, \beta}}
\newcommand{\Lt}{{\cal L}(A,\alpha, \beta)}
\newcommand{\Le}{{\cal L}(E,\alpha, \beta)}
\newcommand{\rM}{{\cal M}_{n}(\mathbb{R})}
\newcommand{\cM}{{\cal M}_{n}(\mathbb{C})}
\newcommand{\M}{{\cal M}_{n}}
\newcommand{\F}{_{_F}}

\title{ The generalized Levinger transformation  }

\author {M. Adam and J. Maroulas$^{1}$}
\maketitle

\vspace{0.2cm}

\footnotetext[1]{Department of Mathematics, National Technical
University, Zografou Campus, Athens 15780, GREECE.
E-mail:maroulas@math.ntua.gr.  Research supported by a grant of
the EPEAEK project Pythagoras II. The project is co-funded by the
European Social Fund $(75^{0}/_0)$ and Greek National Resources
$(25^{0}/_{0}).$
 }
 \vspace{0.2cm}

\begin{abstract}
 In this paper, we present new results relating the numerical range of a matrix $\,A\,$ with
generalized Levinger transformation $\,\Lt=\al H_A+\be S_A,\,$ where
$\,H_A\,$ and $\,S_A, \,$ are, respectively the Hermitian and
skew-hermitian parts of $\,A.\,$ Using these results, we, then
derive expressions for eigenvalues and eigenvectors of the perturbed
matrix $\,A+{\cal L}(E,\al,\be),\,$ for a fixed matrix $\,E\,$ and
$\,\al,\,$ $\be\,$ are real parameters.
\end{abstract}

\vspace{0.15cm}


 \noindent
\section{Introduction}
 \noindent

 Let $\,\cM\, $ $\,(\rM)\,$ be the algebra of all $\,n \times
n\,$ complex (real) matrices, and let $\,A\,\in \,\cM .\,$ The
{\it numerical range} of $\,A,\,$ also known as the {\it field of
values} \cite{HJ}, is the set
$$
NR[A]\,=\,\{\,x^{*}Ax\,\in \, \mathbb{C}\,:\;\;\;x\,\in
\,\mathbb{C}^{n}\,\,\mbox{ with }\;\;x^{*}x\,=\,1\,\}.
$$
The numerical range $\,NR[A]\,$ is a compact and convex subset of
$\,\mathbb{C},\,$ that contains the {\it spectrum} $\,\sigma(A)\,$
of $\,A.\,$ If $\,\lambda \in \sigma (A) \cap \partial NR[A]\,$ with
multiplicity $\,m,\,$  then $\,\lambda \,$ is a {\it normal
eigenvalue}, i.e., $\,A\,$ is unitarily similar to $\,\lambda I_{m}
\oplus B,\,$ with $\,\lambda \notin \sigma(B).\,$ Clearly, when
$\,A\,$ is normal, then its eigenvalues are normal and
$\,NR[A]\,=\,\mbox{Co}\{\,\sigma(A)\,\},\,$ where
$\,\mbox{Co}\{\,\cdot \,\}\,$ denotes the convex hull of the set.
For any $\,A\,\in \,\cM ,\,$ if we write $\,A\,=\,H_A +  S_A,\,$
where
$$\,H_A\,=\,\displaystyle{\frac{A + A^{*}}{\,2\,}}\;\;\;\,\mbox{ and } \;\;\;\,S_A\,=\,\displaystyle{\frac{A - A^{*}}{\,2\,}},\,$$
are the Hermitian and the skew-hermitian parts of $\,A\,$
respectively, then
\begin{eqnarray*}
\label{e1}
 \mbox{Re}\,NR[A] = NR[H_A]\;\;\;\;\mbox{and
} \;\;\;\;  {\bf i}\,\mbox{Im}\,NR[A]= NR[\,S_A].
\end{eqnarray*}
Moreover, if $\,P\,$ is any $\,n \times m\,$ matrix with $\,n\geq
m\,$ and $\,P^{*}P=I_m,\,$ then
\begin{eqnarray*}
\label{e3} \,NR[P^{*}AP]\,\subseteq \, NR[A] ,
\end{eqnarray*}
and the equality holds only for $\,m =n.\,$

Given a matrix $\,A\,\in \,\cM ,\,$ we define  {\it
$^{''}$generalized Levinger transformation} $^{''}$ of $\,A\,$ as
the double parametrized family of matrices
\begin{eqnarray}
\label{e2}
 \Lt\,=\,\al H_A + \be S_A= \frac{\alpha + \beta}{2}A  +\frac{\alpha - \beta}{2}A^{*},\;\;\;\;\mbox{with}\;\;\;\;  \al, \be  \,\in \,
 \mathbb{R}.
\end{eqnarray}
In (\ref{e2}), for $\,\al = 1\,$ and $\,\be = 2t-1,\,$ $\,t\in
\mathbb{R},\,$ we have
$${\cal L}(A,1, 2t-1)=tA+(1-t)A^{*},$$
i.e., $\,{\cal L}(A,1, 2t-1)\,$ is just the ordinary Levinger's
transformation, which has been studied
 in \cite{Fi, L, MPT, McPT}.

Moreover, for $\,\al = 2t-1,\; t\in \mathbb{R}\,$ and $\,\be = 1,\,$
we have
\begin{eqnarray}
\label{newlev} {\cal L}(A,2t-1, 1)\, =\,
tA+(t-1)A^{*}.
\end{eqnarray}
The equation (\ref{newlev}) is a different formulation of Levinger
transformation, where in $\,{\cal L}(A,2t-1, 1)\,$ the {\it
difference} of coefficients $\,t\,$ and $\,t-1\,$ of matrices is
equal to unity (skew convex expression). Note, that
\begin{eqnarray*}
\label{lev5} {\cal L}({\bf i}A,1, 2t-1)= {\bf i}{\cal
L}(A,2t-1,1).
\end{eqnarray*}
Clearly, for every $\,\al, \be \,\in \, \mathbb{R},\,$ we have from
(\ref{e2})
\begin{eqnarray*}
\label{e4} \,H_{\Lt}\,=\,\al H_A,\;\;\;\;\;\,S_{\Lt}\,=\,\be
S_A.\,
\end{eqnarray*}
Hence,
\begin{eqnarray*}
& &\mbox{Re}\,NR[\Lt]=NR[H_{\Lt}]= \al NR[H_A]=\al \, \mbox{Re}\,NR[A]  \label{e5}\\
& &  {\bf i}\,\mbox{Im}\,NR[\Lt]=NR[\,S_{\Lt}]= {\bf i}\be  \,
 \mbox{Im}\,NR[A], \nonumber
\end{eqnarray*}
 and consequently
\begin{eqnarray} \label{e6}
NR[\Lt]\,=\,\{\,\al x+{\bf i}\beta
 y\;\,:\;\;\, x,y\,\in\,\mathbb{R}\,,\;\;\mbox{with }\;\,x+{\bf i}y \in
NR[A]\,\}.
\end{eqnarray}
Moreover, the boundary of $\,NR[\Lt]\,$ is given by
\begin{eqnarray} \label{epar4}
\partial NR[\Lt]\,=\,\{\,\al x+{\bf i}\beta
 y\;\,:\;\;\, x,y\,\in\,\mathbb{R}\,,\;\;\mbox{with }\;\,x+{\bf i}y \in
\partial NR[A]\,\}.
\end{eqnarray}

For $\,A\,\in \,\rM, \,$ since $\,NR[A]=NR[ A^{T}],\,$ $\,NR[\Lt]\,$
is symmetric with respect to the real axis and
\begin{eqnarray*}
NR[\Lt]=NR[ {\cal L}(A^{*},\al, -\be)],
\end{eqnarray*}
the domain of $\,\be\,$ can be reduced to $\,[ 0\,,\,+\infty\,).\,$
Additionally, if $\, 0 < \beta_1 <\beta_2 , \,$ then due to
(\ref{e6}), (\ref{epar4}) and the symmetry of the numerical ranges
with respect to the real axis, we have (in some sense) a vertical
dilation, i.e., for $\,z_1 \in
\partial NR[{\cal L}(A,\al, \be_1) ],\,$ $\,z_2 \in
\partial NR[{\cal L}(A,\al, \be_2)],\,$ holds
$$
|z_1| = \sqrt{ \al^{2}(x^{*}H_Ax)^{2} + \be_1^{2} (x^{*}S_Ax)^{2}\,}\, < \,\sqrt{ \al^{2}(x^{*}H_Ax)^{2} + \be_2^{2} (x^{*}S_Ax)^{2}\,}= |z_2|,\,
$$
and consequently
\begin{eqnarray*}
\label{eqsub} \,NR[ {\cal L}(A,\al, \be_1)] \subset NR[ {\cal
L}(A,\al, \be_2)].\,
\end{eqnarray*}
{\it Example.} For $\,A=\left[%
\begin{array}{ccc}
  1 & 3 & 4 \\
  2 & 7 & -6 \\
  -1 & 3 & 5 \\
\end{array}%
\right],\,$ the numerical ranges of $\,A\,$ and $\,{\cal L}(A,\al,
\beta)\,$ are illustrated in the following figures.  On the right,
due to (\ref{e6}), the vertical dilation is presented only when
$\,a\,$ is fixed (here $a=0.4$) and $\,\beta\,$ is altered ($0.5
\leq \beta \leq 1.2$), reminding that, this property holds only for
the ordinary Levinger's transformation. Otherwise,  $\,NR[\Lt]\,$ is
moved as it is shown on the left figure, for the values $\al=-0.9,
\beta=0.8,\,$ $\al=1.3, \beta=0.6\,$ and $\al=1.4, \beta=1.3.$
\begin{figure}[ht]
\begin{center}
\vspace*{-0.7mm} \mbox{ \hspace*{-5mm}
\epsfig{file=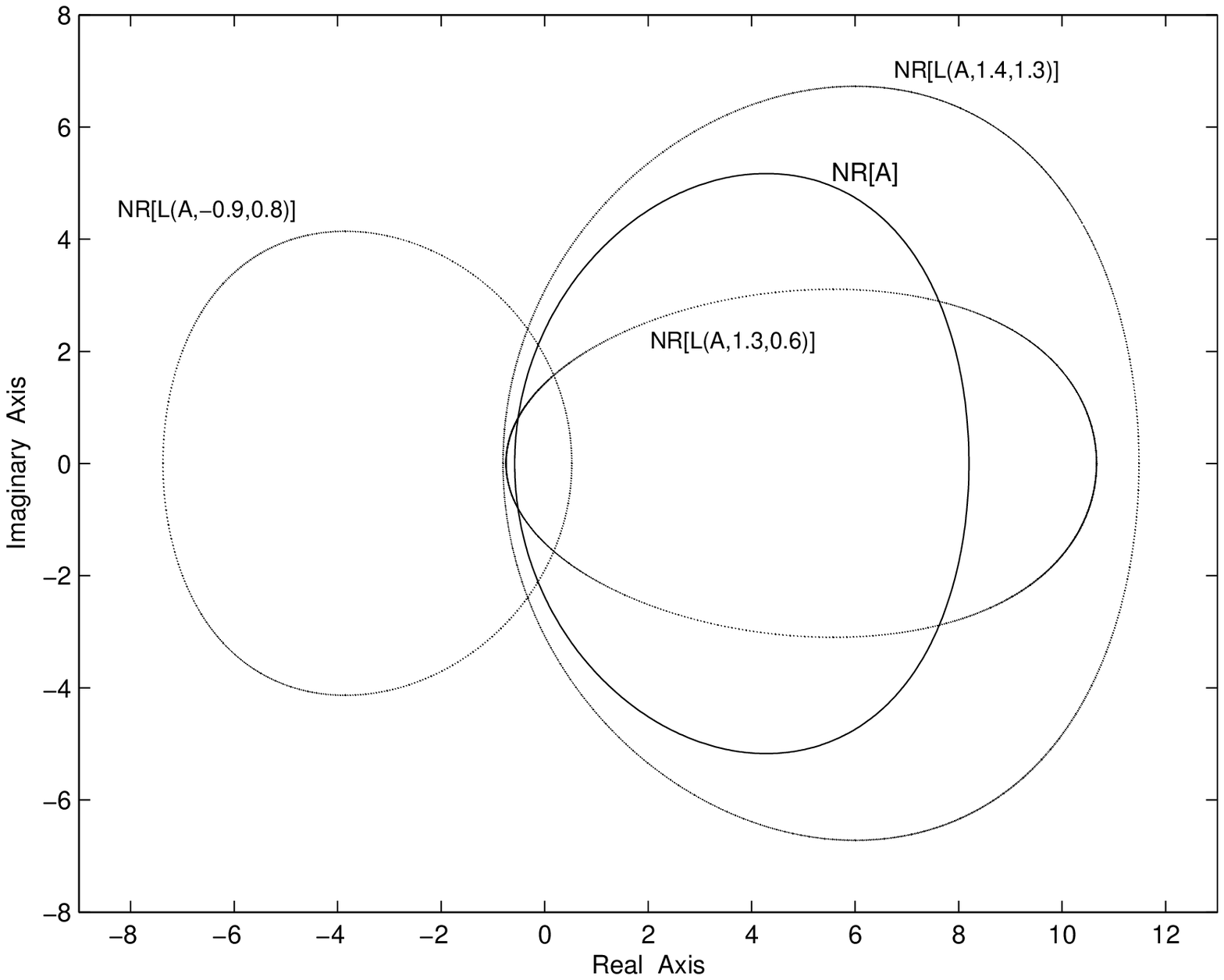,height=68mm,width=72mm,clip=}~~~~~}
\epsfig{file=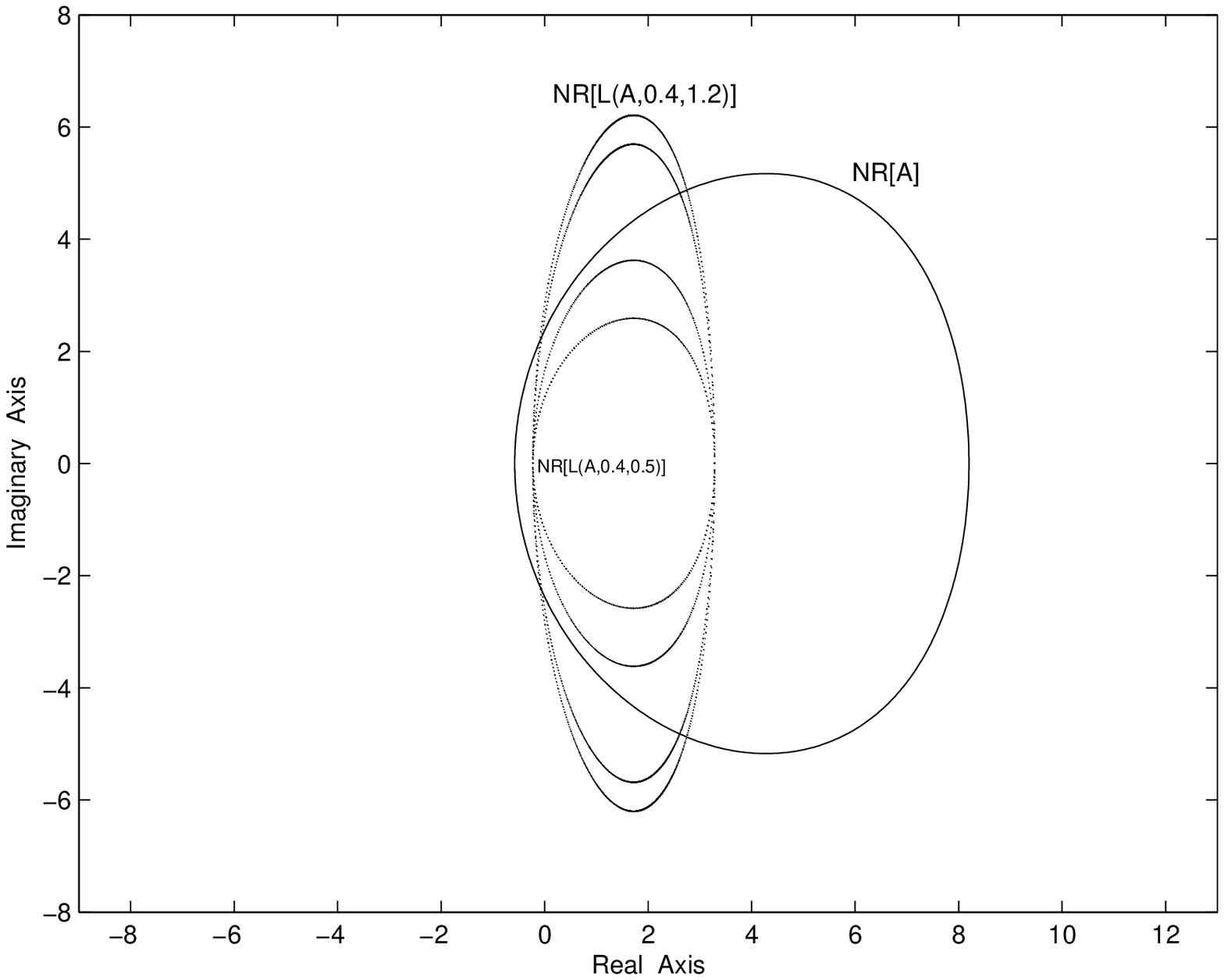,height=68mm,width=72mm,clip=}
\end{center}
\end{figure}
\vspace{0.2cm}
\ \\

The paper is devoted to the study of the generalized Levinger
transformation of a matrix. Specifically, we establish an
interesting relationship between the numerical range of a matrix
$\,A\,$ and its generalized Levinger transformation. This
relationship is then used to obtain results on the eigenvalue and
eigenvector of the perturbed matrix of the form $ \,A+{\cal
L}(E,\al,\be),\,$ where $\,E\,$ is fixed and $\,\al,\,$ and
$\,\beta\,$ are small parameters.

Our motivation for such study comes from the fact that a great deal
of effort has been made in the literature to establish bounds on the
eigenvalues of a perturbed matrix. For results on this topic, see
\cite{X} and the well-known books on linear and numerical linear
algebra by Datta \cite{D}, Stewart and Sun \cite{SS}, Kato \cite{K},
Lancaster and Tismenetsky \cite{LT}, and Bhatia \cite{B}.

This paper is divided in two parts. The first part contains
geometric properties of numerical range of $\,\Lt.\,$ Also bounds
are given for real and imaginary parts of eigenvalues of $\,\Lt.\,$
This provides us a framework to study variation of the spectrum of
$\,\Lt.\,$ In the second part, we use the Levinger transformation
for a fixed matrix $\,E\,$ as a perturbation matrix, whose activity
on a matrix $\,A\,$ depends {\it only} on the parameters $\,\al\,$
and $\,\beta.\,$ First, we formulate a necessary and sufficient
condition for a normal matrix to remain normal, under a perturbation
by a symmetric and rank one matrix.  Next, we present an
approximation of a perturbed eigenpair of a diagonalizable matrix
$A$ using two parameters, which generalizes a known result in
\cite[p. 183]{SS}, where the eigenvector of perturbed eigenvalue is
not mentioned and even the perturbed eigenvector is investigated in
\cite[p. 431]{D}, without giving further details for the
perturbation of the corresponding eigenvalue. Further we simplify
these formulae using the notion of generalized inverse extending the
corresponding result in \cite{Fi}. As an application, we present a
sufficient condition such that the perturbed eigenpair is first
order approximation of the corresponding simple eigenpair of initial
matrix $\,A,\,$ and we give two numerical examples to illustrate our
results.

 \noindent
\section{Geometric properties}

\begin{prop}
\label{line} ~{\bf a.}~Let $\,A\,\in \cM.\,$ The image of a line
segment $\,\epsilon_{A} \in NR[A]\,$ by the Levinger transformation,
is a line segment $\,\epsilon_{\cal L} \in NR[\Lt].\,$
\newline ~{\bf b.}~If $\,A\,\in \rM\,$ and $\,NR[A]\,$ is an
ellipse, then $\,NR[\Lt],\,$ for $\,\al \neq 0,\,$ is also an ellipse.
\end{prop}
\pr~~{\bf a.}~ For any $\,x_1+{\bf i}y_1,x_2+{\bf i}y_2 \in
NR[A],\,$ and $\, t \in [0,1],\,$ we observe that
$$
\al [ (1-t) x_1 + t x_2] + \beta [ (1-t){\bf i} y_1 + t{\bf i}
y_2] = (1-t)(\al x_1+{\bf i} \beta y_1) + t (\al x_2+{\bf i} \beta
y_2),
$$
where $\,\al x_1+{\bf i}\beta y_1\,$ and $\,\al x_2+{\bf i}\beta
y_2\,$ lie in the convex set $\,NR[\Lt].\,$ Hence,  the proof of
{\bf a} follows readily.
\newline
~{\bf b.}~Consider that $\,NR[A]\,=\, \{\,x + {\bf i} y
\,:\,\frac{\,x^{2}\,}{\,c^{2}\,} \,+\,\frac{\,y^{2}\,}{\,k^{2}\,}
\,\leq \,1,\, \;\mbox{ with }\, c> k>0\,\,\}.\,$ For $\, \al \neq
0,\,$ and $\, \be > 0, \,$ changing the variables $\,x=\al^{-1}
X,\,$ $\,y=\be^{-1}Y,\,$ by (\ref{e6}) we have that the boundary
$\,\partial NR[\Lt]\,$ is the ellipse
$\,\frac{\,X^{2}\,}{\,\al^{2} c^{2}\,}+
\frac{\,Y^{2}\,}{\be^{2}k^{2}}\,=\,1 . \,$ The foci are on the
real axis, when $\,\al \in \mathbb{R}\backslash [\,-\frac{\,\be
k\,}{\,c\,}\,,\,\frac{\,\be k\,}{\,c\,}\,],\,$ otherwise, they lie
on the imaginary axis.
 \qed

\begin{prop}
\label{pextra} ~Let $\,A\in \cM,\,$ $\,B \in {\cal
M}_{m}(\mathbb{C}),\,$ and let $\,NR[B]\,$ be a polygon
circumscribed to $\,NR[A].\,$  Then the geometric relationship of
numerical ranges $\,NR[\Lt]\,$ and $\,NR[{\cal L}(B,\al, \be)]\,$
remains the same, for  $\, \al,\be \in \mathbb{R}- \{0\}.\,$
\end{prop}
\pr~By Proposition \ref{line}{\bf a} and (\ref{epar4}), $NR[{\cal
L}(B,\al,\be)]$ is a convex polygon. Moreover, it is easy to see
that there do not exist other common boundary points of
$\,NR[\Lt]\,$ and $\, NR[{\cal L}(B,\al, \be)],\,$ except those,
which correspond to $\,z \in \partial NR[A] \cap
\partial NR[B].\,$ \qed

\begin{prop} \label{p1} ~Let $\,A\,\in \cM,\,$
and $\, \al,\be \in \mathbb{R}- \{0\}.\,$
\begin{description}
  \item[a.] If $\,A\,$ is a normal matrix, then $\,\partial NR[\Lt]\,$
  is a $k-$polygon, as $\,\partial NR[A].\,$
  \item[b.] $\,NR[\Lt]\,\cap \mathbb{R}\,=
  \,\mathbb{R}\,\cap\,NR[\al A].\,$
\end{description}
\end{prop}
\pr ~{\bf a.} The proof of this part follows readily by Proposition
\ref{line}{\bf a} and from the observation that $\,A\,$ is unitarily
similar to $\,diag\{\,x_1+{\bf i}y_1, \,x_2+{\bf i}y_2, \ldots ,
x_n+{\bf i}y_n \}\,$  if and only if $\,\Lt\,$ is unitarily similar
to $\,diag\{\,\al x_1+{\bf i}\beta y_1, \,\al x_2+{\bf i}\beta y_2,
\ldots , \al x_n+{\bf i}\beta y_n \}.$
\newline
{\bf b.}~Clearly, if $\,\al x^{*}Ax \in NR[\al A] \cap
\mathbb{R},\,$ then $\,x^{*}\Lt x = \frac{\, \al + \be
\,}{\,2\,}\,x^{*}A x + \frac{\, \al - \be \,}{\,2\,}\, x^{*}A^{*}x =
\al x^{*}Ax, \,$ concluding that $\,NR[\al A] \cap \mathbb{R}
\subset NR[\Lt] \cap \mathbb{R}.\,$ Moreover, for $\,\al,\,\beta \in
\mathbb{R} \backslash \{ 0\},\,$ because
\begin{eqnarray*} \label{e8n}
\al A=\frac{\, \al + \be \,}{\,2 \be \,}\,\Lt - \frac{\, \al - \be
\,}{\,2 \be \,}\,{\cal L}^{*}(A,\al,\be) ,
\end{eqnarray*}
in an analogous way, we obtain $\,NR[\Lt] \cap \mathbb{R} \subset
NR[\al A] \cap \mathbb{R}.\,$
 \qed \ \\
\ \\
{\bf Remark }
\newline ~ The eigenvalue $\,\la \in \sigma(A)\,$ is {\it normal} if and only if
$\,\la_{{\cal L}}\,$ is a normal eigenvalue of $\,\Lt.\,$ \newline
In fact, from (\ref{e2}) and the relationship $ U^{*}A U = \lambda
I_{m} \oplus B,\, $ where $\, \lambda \notin \sigma(B),\,$ it
follows that $\,U^{*} \Lt U = \la_{{\cal L}}I_{m} \oplus {\cal
L}(B,\al,\be) .\, $
 \\

 In the following proposition we present a compression of $\,NR[\Lt],\,$  when $\,A\,$ is a normal matrix, based on our results in \cite{AM}.
\begin{prop}
\label{pcompress} ~Let $\,A\,\in \cM\,$ be a normal matrix and let
the polygon
$\,<\la_{(A)_{1}}\,,\,\la_{(A)_{2}}\,,\,\ldots\,,\,\la_{(A)_{k}}>\,$
be the numerical range of $\,A.\,$ If $\,x_j\,$ is a corresponding
eigenvector of $\,\la_{(A)_{j}}, \;j = 1,2,\,\ldots, k,\,$ and
$\,\upsilon \,=\, \sum_{ j =1 }^{k} \upsilon_j x_j\,$ is a unit
vector, denoting by $\,E\,=\,span\{\,\upsilon \,\}\,$ and
$\,E_{W}^{\bot}\,$ the orthogonal complement of $\,E\,$ with respect
of $\,W\,=\,span \{\,x_1\,,\,x_2\,,\,\ldots\,,\,x_{k}\,\},\,$  then
$$
NR[P^{*}\Lt P]= NR [{\cal L}(P^{*}AP, \al, \be)]\,\subset \,<
\la_{({\cal L})_{1}}\,,\,\la_{({\cal
L})_{2}}\,,\,\ldots\,,\,\la_{({\cal L})_{k}}
>,\,$$
where $\,
P\,=\,[\,w_1\,\,\,\,w_2\,\,\,\,\ldots\,\,\,\,w_{k-1}\,],\,$ and
$\,w_1\,,\,w_2\,,\,\ldots\,,\,w_{k-1}\,$ is an orthonormal basis of
$\,E^{\bot}_{W}.\,$ Moreover, $\,\partial NR[P^{*}\Lt P]\,$ is
tangential to the edges of the polygon $\,<\la_{({\cal
L})_{1}}\,,\,\la_{({\cal L})_{2}}\,,\,\ldots\,,\,\la_{({\cal
L})_{k}}>\,$ at the points
\begin{eqnarray}
\label{epoint} \mu_{ ({\cal L}) \tau } \,=\,\al
\mbox{Re}\,\mu_{(A)\tau} + i\,\be \mbox{Im}\, \mu_{(A)\tau},\;
\;\;\;(\tau\,=\,1, \ldots,k)
\end{eqnarray}
where
\begin{eqnarray*}
\mu_{(A)\tau}=\frac{|\upsilon_{\tau +
1}|^{2}\,\,\la_{(A)\tau}\,+\,|\upsilon_{\tau}|^{2}\,\,\la_{(A)_{\tau
+1}}}{|\upsilon_{\tau +1}|^{2} + |\upsilon_{\tau}|^{2}}
\;\;\;\;\;(\tau = 1, \ldots ,k-1),\;\;\;\;\mu_{(A)_{k}} \,=\,
\frac{|\upsilon_1|^{2}\,\,\la_{(A)_{k}}\,+|\upsilon_{k}|^{2}\,\,\la_{(A)_{1}}}{|\upsilon_1|^{2}
+ |\upsilon_{k}|^{2}}.
\end{eqnarray*}
\end{prop}
\pr ~By the equation (\ref{e2}) we have $ {\cal
L}(P^{*}AP,\al,\be)\, =\,P^{*}\Lt P,$ and it is evident the
equality of the numerical ranges. Moreover, $\,{\cal
L}(A,\al,\be)\,$ is normal, and it is known in \cite{AM} that, $\,
\partial NR[P^{*}\Lt P]\,$ tangents to the edges of the
polygon $\,<\la_{({\cal L})_{1}}\,,\,\la_{({\cal L})_{2}}\,,\ldots
,\,\la_{({\cal L})_{k}}>$ at the points
$$\mu_{({\cal L})\tau}=\frac{|\upsilon_{\tau +
1}|^{2}\,\,\la_{({\cal
L})\tau}\,+\,|\upsilon_{\tau}|^{2}\,\,\la_{({\cal L})_{\tau
+1}}}{|\upsilon_{\tau +1}|^{2} + |\upsilon_{\tau}|^{2}}
\;\;(\tau\,= 1, \ldots ,k-1)\;\;\;;\;\;\; \mu_{({\cal L})_{k}}
\,=\, \frac{|\upsilon_1|^{2}\,\,\la_{({\cal
L})_{k}}\,+|\upsilon_{k}|^{2}\,\,\la_{({\cal
L})_{1}}}{|\upsilon_1|^{2} + |\upsilon_{k}|^{2}}.$$ Since, the
eigenvalues of $\,A\,$ and $\,\Lt\,$ are related by $\,\la_{({\cal
L})\tau}= \al\,\mbox{Re} \lambda_{(A)_{\tau}}+i \be\,\mbox{Im}
\lambda_{(A)_{\tau}},\,$  the equation (\ref{epoint}) is verified.
\qed \ \\

In the following applying the results of Rojo and Soto in
\cite{RS} to $\Lt, \,$ one obtains bounds for the real and the
imaginary parts of the eigenvalues of $\,\Lt.\,$
\begin{them}
\label{ineq} ~Let the matrix $\,A\in \rM\,$ and $\,\la_j \in \sigma(\Lt).\,$ Then for each $\,\la_j \,$ we have
\begin{eqnarray}
\label{inre} \left|\mbox{Re} \la_j - \frac{\al tr
(H_A)}{\,n\,}\right| \, \leq \, |\al| \left[ \frac{n-1}{\,n\,}
\left( \|H_A\|^{2}_F -\, \frac{ \be^{2}\|S_AH_A-H_AS_A\|^{2}_F}
{\,3\left(\al^{2}\|H_A\|^{2}_F +\be^{2}\|S_A\|^{2}_F
\right)\,}-\frac{\left[tr(H_A)\right]^{2}}{\,n\,} \right)\,
\right]^{1/2},
\end{eqnarray}
and
\begin{eqnarray}
\label{inim}
\left| \mbox{Im} \la_j \right| \, \leq \, |\be|\, \left [ \frac{n-1}{\,n\,} \left( \|S_A\|^{2}_F - \, \frac{ \al^{2}\|S_AH_A-H_AS_A\|^{2}_F}
{3\left(\al^{2}\|H_A\|^{2}_F +\be^{2}\|S_A\|^{2}_F\right)\,} \,  \right) \right]^{1/2},
\end{eqnarray}
where $\,\|\cdot\|_F\,$ denotes the Frobenius norm.
\end{them}
\pr~ Observe that $\, tr(\Lt) = tr(\al H_A+\be S_A)=\al tr(H_A),\,$ and
\begin{eqnarray}
\label{rs}
& & \|\Lt{\cal L}^{T}(A, \al, \be )-  {\cal L}^{T}(A, \al, \be )\Lt\|_F ~~~~~~~~~~~~~~~~~~~~~~~~~~~~~~ \nonumber  \\
& & ~~~~~~~~~~ = \|(\al H_A+\be S_A)(\al H_A^{T} + \be S_A^{T})-(\al H_A^{T} + \be S_A^{T})(\al H_A+\be S_A)\|_F \nonumber \\
&& ~~~~~~~~~~ =2|\al \be| \|S_AH_A-H_AS_A\|_F.
\end{eqnarray}
Since
$$
tr(H_AS_A^{T})=-tr(H_AS_A)=-tr\left(\frac{(A+A^{T})(A-A^{T})}{\,4\,}\right)=-\frac{\,1\,}{\,4\,}\left(
tr(A^{2})-tr\left[(A^{T})^{2}\right]\right)=0
$$
we have
\begin{eqnarray}
\label{rs1}
\|\Lt\|_F^{2}& = & tr[ \Lt {\cal L}^{T}(A, \al, \be ) ] = tr \left[ (\al H_A + \be S_A)(\al H_A^{T}+ \be S_A^{T}) \right] \nonumber \\
& = &
\al^{2}\|H_A\|^{2}_F + \be^{2}\|S_A\|^{2}_F
\end{eqnarray}
and
\begin{eqnarray}
\label{rs2}
tr[{\cal L}^{2}(A, \al, \be )] = tr \left[ (\al H_A + \be S_A)(\al H_A + \be S_A) \right]=\al^{2}\|H_A\|^{2}_F - \be^{2}\|S_A\|^{2}_F,
\end{eqnarray}
and consequently by (\ref{rs1}) and (\ref{rs2})
\begin{eqnarray}
\label{equa}
& & \|\Lt\|^{2}_F + tr({\cal L}^{2}(A, \al, \be ))=2\al^{2}\|H_A\|^{2}_F, \nonumber \\
\\
& & \|\Lt\|^{2}_F - tr({\cal L}^{2}(A, \al, \be ))=2\be^{2}\|S_A\|^{2}_F.\nonumber
\end{eqnarray}
Therefore, if we substitute (\ref{rs}),  (\ref{rs1}), and
(\ref{equa}) in the relationships of Theorem 7 in \cite{RS}, we have
$$
\left|\mbox{Re} \la_j - \frac{ tr (\Lt)}{n}\right| \leq \sqrt{ \frac{n-1}{n} \left( \frac{\|\Lt\|^{2}_F + tr({\cal L}^{2}(A,\al,\be ))}{2}
-\frac{ \nu(\Lt)^{2} }
{12 \| \Lt \|^{2}_F }-\frac{[tr(\Lt)]^{2}}{n} \right ) }
$$
and
\begin{eqnarray*}
\left| \mbox{Im} \la_j \right|  \leq  \sqrt{ \frac{n-1}{\,2n\,} \left( \|\Lt\|^{2}_F - tr({\cal L}^{2}(A, \al, \be ))-
\frac{ \nu(\Lt)^{2}}
{6 \| \Lt \|^{2}_F\,}  \right) }\,,
\end{eqnarray*}
where $\,\nu(\Lt)= \|\Lt{\cal L}^{T}(A, \al, \be )-{\cal L}^{T}(A,
\al, \be )\Lt\|_F,\,$ thus we obtain the bounds for the real and the
imaginary part for each eigenvalue of $\Lt$ in (\ref{inre}) and
(\ref{inim}). \qed
\newpage
\noindent
\section{Application to perturbation theory}
\noindent

 The question $^{''}${\it how close is a matrix} $\,M\,$ {\it to
being normal }$^{''}$, it is known that it is evaluated by the {\it
normality distance} $\,\|AA^{T}-A^{T}A\|_p\,$ of $\,p\,$ norm. Since
various matrix norms are equivalent, using the Frobenius norm we
have:
\begin{prop}
\label{normalsym} ~Let $\,N\in \rM\,$ be a normal matrix and for a
nonzero vector $\, x \in \mathbb{R}^{n},\,$ let $\,E = xx^{T}.\,$ If
$\,x\,$ is not eigenvector of $\,N\,$ corresponding to a real
eigenvalue, the matrix $\,M=N+E\,$ is normal if and only if $\,N\,$
is symmetric.
\end{prop}
\pr~ For the symmetric matrix $\,E= xx^{T},\,$ clearly
$\,E^{2}=\|x\|^{2}E,\,$ and the normality distance of $\,M\,$ is
equal to
\begin{eqnarray*}
\|MM^{T}-M^{T}M\|\F & = & \|(N+E)(N^{T}+E)-(N^{T}+E)(N+E)\|\F \\
& =& \|NE+EN^{T}-N^{T}E-EN \|\F = \|R+R^{T}\|\F,
\end{eqnarray*}
where $\,R=NE-EN.\,$ Since
\begin{eqnarray*}
tr(R^{2})& = & 2 \left[\, tr(NENE)-tr(EN^{2}E) \right]= 2 \left[
(x^{T}Nx)^{2}-\|x\|^{2} (x^{T}N^{2}x) \right] \\
tr(RR^{T})& = & 2 \left[ \,\|x\|^{2}tr(NEN^{T})-tr(ENEN^{T})
\right]= 2 \left[ \|x\|^{2}(x^{T}N^{T}Nx)- (x^{T}Nx)^{2} \right]
\end{eqnarray*}
and $\,tr\left[(R^{T})^{2}\right] = tr(R^{2}),\,$ we have :
$$
\|R+R^{T}\|^{2}\F = tr\left[( R+R^{T})^{2}\right]=
tr(R^{2})+tr\left[(R^{T})^{2}\,\right]+2tr(RR^{T})=
4\|x\|^{2}\,x^{T}(N^{T}N-N^{2})x.
$$
Hence, the matrix $\,M\,$ is normal if and only if
$\,N^{T}N=N^{2}.\,$ This equation is equivalent to
\begin{eqnarray}
\label{normality}
 \overline{D}D=D^{2}, \end{eqnarray}
where $\,D\,$ is diagonal and unitary similar to $\,N,\,$ i.e.,
$\,N=UDU^{*}.\,$ Thus, by (\ref{normality}), $\,D\,$ is real and
$\,N\,$ is symmetric, since it is unitary similar to a real
diagonal matrix. \qed
\vspace{0.4cm}

\begin{coro}
~ Let $\,N=diag(N_1,\,N_2,\,\ldots,\,N_{\tau})\in \rM\,$ be a normal
matrix and for a nonzero vector $\,x=\left[%
\begin{array}{ccc}
  x_1^{T} & \ldots & x_{\tau}^{T} \\
\end{array}%
\right] \in \mathbb{R}^{n},\,$ with all $\,x_j \neq 0,\, $ let
$\,E = xx^{T}.\,$ If $\,x_j \, $ is not eigenvector of $\,N_j,\,$
$\,(j=1, \ldots ,\tau),\,$ corresponding to a real eigenvalue, the
matrix $\,M=N+E\,$ is normal if and only if $\,N\,$ is symmetric.
\end{coro}

It is worth notice that, the result in Proposition \ref{normalsym}
is combined by the special form of $\,E,\,$ but it is interesting to
look at more general perturbations, investigating how main
properties of $\,N\,$ are influenced. For this, we consider $\,E \in
\rM , \,$ since
$$
\|MM^{T}-M^{T}M\|_F = \|2H_{[N,E^{T}]} +EE^{T}-E^{T}E\|_F,
$$
where $\,[N,E^{T}] = NE^{T}-E^{T}N,\,$ we conclude that, the
normality distance of $\,M\,$ is related to the normality distance
of $\,E.\,$  Hence, an outlet is to investigate {\it if some
properties of perturbed normal matrices remain}.

Let the matrix $\,A \in \cM\,$ be {\it diagonalizable} (keeping
the property of $\,N\,$) and $\,E \in \rM\,$ be fixed, without
giving any attention to $\,\|E\|.\,$ Consider the matrix
\begin{eqnarray}
\label{eqM} M_{\al,\be} = A + {\cal L}(E,\al,\be)= A + \al H_E +
\be S_E ,
\end{eqnarray}
where $\,\al,\be \in \mathbb{R}\,$ are small enough varying
parameters. Clearly in (\ref{eqM}), $\, M_{\al,\be}$ is continuous
differentiable and the Hermitian and the skew-hermitian parts of
$\,E\,$ influence independently the matrix $\,A.\,$ Especially, when
$\,A\,$ is normal, $\,H_E\,$ and $\,S_E\,$ alter $\,H_A\,$ and
$\,S_A\,$ separately.

Denote by $\,\lab\,$ an eigenvalue of $\,M_{\al,\be} \,$ in
(\ref{eqM}) and by $\,\uab\,$ and $\,\oab,\,$ the corresponding
right and left eigenvectors, i.e., $\,(M_{\al,\be}- \lab I) \uab
=0,\,$ $\,\oab^{*}(M_{\al,\be}- \lab I) =0.\,$ Since the
coefficients of characteristic polynomial $\,det(\la I -
M_{\al,\be}) \,$ are polynomials of two variables $\,\alpha,
\beta\,$ and $\,\lab\,$ is continuous function of these
coefficients, for $\,\al = \be = 0\,$ the perturbed eigenvalue
$\,\lab \,$ is equal to a semisimple eigenvalue $\,\la_i\,$ of
$\,A,\,$ and the eigenvectors are : $\,\uab = \upsilon_i , \,$
$\,\oab = \omega_i, \,$ where $\,\upsilon_i\,$ and $\,\omega_i\,$
are the right and left eigenvectors of $\,\la_i \,$ for the matrix
$\,A.\,$ We remind the readers that an eigenvalue is called {\it
semisimple}, when it is a simple root of the minimal polynomial of
matrix. Moreover, $\,\lab\,$ and $\, \uab,\,\oab \,$ are continuous
functions of $\,\alpha, \beta \,$ and partial differentiable, but
might have rather singularities on total differentiability \cite[p.
116]{K}. For further details we refer to \cite{K} and \cite[Ch.
11]{LT}. We will now give a result on the sensitivity of eigenvalues
and eigenvectors of perturbed matrix $\,M_{\al,\be}\,$ in
(\ref{eqM}) in the neighborhood of $\,\la_i\,$ in relation with the
remaining eigenvalues and eigenvectors.

\begin{them}
\label{normaleig} ~Let the matrix $\,A \in \cM\,$ be diagonalizable
and let $\,\upsilon_j\,$ and $\,\omega_j\,$ be the right and left
eigenvectors of $ \,A\,$ corresponding to $\,\la_j \in \sigma(A).\,$
If the eigenpair $\,\left(\lab\,,\,\uab \right)\,$ has continuous
second order partial derivatives in the neighborhood of $\,\la_i\,$
and $\,\upsilon_i$, then :
\begin{eqnarray}
\label{eigen1} \lab = \la_i + \frac{ \omega^{*}_i \Le \upsilon_i
}{ s_i } + \sum_{ k \neq i} \frac{ (\, \omega^{*}_i\, \Le\,
\upsilon_k\,)( \omega^{*}_k \,\Le \,\upsilon_i\,) }{ (\lambda_i -
\lambda_k)\, s_i s_k} + \mathcal{O}\left( \al^{3},\be^{3}\right)~~~~~~~
\end{eqnarray}
\begin{eqnarray}
\label{eigen2} \uab = \upsilon_i + \sum_{k \neq i} \frac{
\,\upsilon_k \, \omega^{*}_k \, \Le \,\upsilon_i }{(\lambda_i -
\lambda_k)\,s_k
 } + \sum_{j \neq i}
\sum_{k \neq i }\frac{ ( \omega_j^{*}\, \Le \, \upsilon_k )(
\omega_k^{*}\,\Le \, \upsilon_i ) }{ ( \lambda_i - \lambda_k) (
\lambda_i -
\lambda_j )\,s_k s_j }\,\upsilon_j ~~~~~~~~  \nonumber \\
- \sum_{j \neq i }\frac{ ( \omega_j^{*}\, \Le \, \upsilon_i )( \omega_i^{*}\,\Le \,
\upsilon_i ) }{ ( \lambda_i - \lambda_j )^{2}\,s_i s_j }\,\upsilon_j +
\mathcal{O}\left( \al^{3},\be^{3}\right), ~~~~~~~~~~~~~~~~~~~ ~~~~~~
\end{eqnarray}
where $\,s_{\ell}=\omega^{*}_{\ell} \,\upsilon_{\ell}.\,$
\end{them}
\pr~ The partial derivatives of the equation $ \left( M_{\al,\be}
- \lab I \right) \uab = 0,\, $ with respect to $\,\al,\,\be,\,$
are
\begin{eqnarray}
 \left( H_E \,- \, \frac{\partial \lab  }{ \partial \al } I \,\right) \uab
+ \left( M_{\al,\be} \, - \,
  \lab I \,\right)
\frac{\partial \uab  }{ \partial \al } = 0 \nonumber \\
\label{deriv1}
\\
\left( S_E \,- \frac{\partial \lab }{
\partial \be } I\, \right) \uab +
\left( M_{\al,\be} - \lab I\, \right) \, \frac{\partial \uab }{
\partial \be } = 0 \nonumber
\end{eqnarray}
Multiplying these by $\,\oab^{*}, \,$ since
$\,\oab^{*}M_{\al,\be}= \lab \oab^{*}$, we have
\begin{eqnarray*}
\oab^{*} \left( H_E \,-\,\frac{\partial \lab }{
\partial \al }\,I \, \right) \uab  = 0, \;\;\;\;\;\;\; \oab^{*} \left( S_E  \,- \frac{\partial \lab  }{
\partial \be } \,I \, \right) \uab = 0.
\end{eqnarray*}
For $\, (\al, \be) \rightarrow ( 0, 0) \,$ the  expressions given
above
\begin{eqnarray}
\label{elabda1} \frac{\partial \lambda_{(\al ,\be)=0}}{
\partial \al } =\lim_{(\al,\be) \rightarrow ( 0, 0)} \frac{\partial \lambda_{\al ,\be}}{
\partial \al }=   \frac{ \omega^{*}_i H_E \upsilon_i  }{ \omega^{*}_i \upsilon_i
},\;\;\;\;\;\;\;\frac{\partial \lambda_{(\al ,\be)= 0} }{
\partial \be } = \lim_{(\al,\be) \rightarrow ( 0, 0)}\frac{\partial \lambda_{\al ,\be}
}{
\partial \be }=    \frac{ \omega^{*}_i S_E \upsilon_i  }{ \omega^{*}_i \upsilon_i
}, \end{eqnarray} and then, the first differential
$\,d\lambda_{\al ,\be}\,$ is equal to
\begin{eqnarray*}
d\lambda_{\al ,\be}\,=\, \al \,\frac{\partial \lambda_{(\al
,\be) = 0}}{
\partial \al } + \be\,\frac{\partial \lambda_{(\al ,\be) = 0}}{
\partial \be } = \al \, \frac{ \omega^{*}_i H_E \upsilon_i }{
\omega^{*}_i \upsilon_i } + \be \, \frac{ \omega^{*}_i S_E
\upsilon_i }{ \omega^{*}_i \upsilon_i}\,=\,\frac{ \omega^{*}_i \Le
\upsilon_i }{ \omega^{*}_i \upsilon_i}.
\end{eqnarray*}
Moreover, the first equality in (\ref{deriv1}) for $\, (\al, \be)
\rightarrow ( 0, 0)\,$ gives
\begin{eqnarray*}
 \left( H_E \,- \, \frac{\partial \lambda_{(\al
,\be) = 0}  }{ \partial \al } I \,\right) \upsilon_i + \left( A \, -
\,
  \lambda_i I \right)
\frac{\partial \upsilon_{(\al ,\be) = 0}  }{ \partial \al } = 0.
\end{eqnarray*}
Since $\,A\,$ is diagonalizable, we can write $\,\frac{\partial
\upsilon_{(\al ,\be) = 0}}{
\partial \al } = \sum_{k=1}^{n} c_k \upsilon_k,\,$
and so the last equality can be written as
\begin{eqnarray*}
 \left( H_E \,- \, \frac{\partial \lambda_{(\al
,\be) = 0}  }{ \partial \al } I \,\right) \upsilon_i + \sum _{ k\neq
i} c_k \left( \lambda_k  -
  \lambda_i  \right) \upsilon_k = 0.
\end{eqnarray*}
Furthermore, multiplying  the above equality by the left eigenvector
$\,\omega_k \,$ of $\,A,\,$ and using the orthogonality of
$\,\omega_k\,$ and $\,\upsilon_i\,$  $\,( k \neq i),\,$ we have
$$\,c_k = \frac{ \omega^{*}_k H_E\upsilon_i\,}{\,
(\lambda_i-\lambda_k) \, \omega^{*}_k \, \upsilon_k },\,\;\;\mbox{
for }\;\;\,k \neq i,\, $$ and consequently,
\begin{eqnarray}
\label{equa1}
\frac{\partial \upsilon_{(\al ,\be) = 0} }{
\partial \al }= \sum_{k \neq i} \frac{ \omega^{*}_k H_E \upsilon_i }{(\lambda_i - \lambda_k)\,\omega^{*}_k \, \upsilon_k
 }\,\upsilon_k.
\end{eqnarray}
Similarly, by the second equality in (\ref{deriv1}), we obtain $\,
(\lambda_i-\lambda_k ) \omega^{*}_k \, \displaystyle{ \frac{\partial
\upsilon_{(\al ,\be) = 0} }{
\partial \be } }\,= \, \omega^{*}_k S_E\upsilon_i  ,\,
$ and thus
\begin{eqnarray}
\label{equa2}
 \frac{\partial \upsilon_{(\al ,\be) = 0} }{
\partial \be }= \sum_{k \neq i} \frac{ \omega^{*}_k S_E \upsilon_i }{(\lambda_i - \lambda_k)\,\omega^{*}_k \, \upsilon_k\,
 }\,\upsilon_k.
\end{eqnarray}
Hence, the differential $\,d\upsilon_{\al ,\be}\,$ can be computed
as
\begin{eqnarray*}
d\upsilon_{\al ,\be} & = & \al \,\frac{\partial \upsilon_{(\al
,\be) = 0}}{
\partial \al } + \be\,\frac{\partial \upsilon_{(\al ,\be) = 0}  }{
\partial \be } \\
& =&  \al \sum_{k \neq i} \frac{ \omega^{*}_k H_E \upsilon_i }{(\lambda_i
- \lambda_k)\,\omega^{*}_k \, \upsilon_k
 }\,\upsilon_k+ \be \sum_{k \neq i} \frac{ \omega^{*}_k S_E \upsilon_i
}{(\lambda_i - \lambda_k)\,\omega^{*}_k \, \upsilon_k
 }\,\upsilon_k \, =   \sum_{k \neq i} \frac{
\,\upsilon_k \, \omega^{*}_k \, \Le \,\upsilon_i }{(\lambda_i -
\lambda_k)\,\omega^{*}_k \, \upsilon_k }.
\end{eqnarray*}
Now, the partial derivatives of the equations in (\ref{deriv1}) with
respect to $\,\al,\be, \,$ are
\begin{eqnarray}
\label{deriv2} 2\left( H_E \, -\, \frac{\partial \lab }{
\partial \al}\,I \right)  \frac{\partial \uab }{
\partial \al}  + \left( M_{\al,\be} -  \lab \,I \, \right)  \frac{\partial^{2} \uab }{
\partial \al^{2}} - \frac{\partial^{2} \lab }{
\partial \al^{2}} \,\uab =0 ~~~~~~~~~~~~~~~~~\nonumber  \\
2\left( S_E \, -\, \frac{\partial \lab }{
\partial \be} \,I \,\right)  \frac{\partial \uab }{
\partial \be} + \left( \, M_{\al,\be}- \lab I\, \right)  \frac{\partial^{2} \uab }{
\partial \be^{2}} - \frac{\partial^{2} \lab }{
\partial \be^{2}} \,\uab =0  ~~~~~~~~~~~~~~~~~ \\
\left( M_{\al,\be} - \lab \,I \right)\frac{\partial^{2} \uab }{
\partial \al \, \partial \be } + \left( H_E -    \frac{\partial \lab }{ \partial \al } \,I  \right)\frac{\partial \uab }{ \partial \be }
+ \left( S_E - \frac{\partial \lab }{
\partial \be} \,I \right) \frac{\partial \uab }{
\partial \al } - \frac{\partial^{2} \lab }{
\partial \al \, \partial \be } \,\uab = 0.\nonumber
\end{eqnarray}
\ \\
Multiplying these expressions by $\,\oab^{*}\,$ and substituting
$\,\frac{\partial \upsilon_{(\al ,\be) = 0} }{
\partial \al},\,\frac{\partial \upsilon_{(\al ,\be) = 0} }{
\partial \be} \,$  from (\ref{equa1}) and (\ref{equa2}),  for $\, (\al,
\be) \rightarrow ( 0, 0), \,$ and noting that $\,\omega^{*}_i
\upsilon_k = 0,\,$ we obtain
\begin{eqnarray}
\label{eigall}
&& \frac{\partial^{2} \lambda_{(\al ,\be) = 0} }{
\partial \al^{2}}=  \frac{\,2\,}{ \omega_{i}^{*} \upsilon_i } \, \left(  \omega_i^{*} H_E  - \omega_i^{*} \frac{\partial  \lambda_{(\al ,\be) = 0}  }{
\partial \al} \right) \frac{\partial  \upsilon_{(\al ,\be) = 0}  }{
\partial \al} =  \frac{\,2\,}{ \omega_{i}^{*} \upsilon_i } \sum_{ k \neq i} \frac{ (\, \omega^{*}_k H_E\, \upsilon_i\,)( \omega^{*}_i H_E \,\upsilon_k\,) }{
(\lambda_i - \lambda_k)\,\omega^{*}_k \upsilon_k }   \nonumber \\
&& \frac{\partial^{2} \lambda_{(\al ,\be) = 0}}{
\partial \be^{2}}=  \frac{\,2\,}{ \omega_{i}^{*} \upsilon_i }\, \left(  \omega_i^{*} S_E - \omega_i^{*} \frac{\partial  \lambda_{(\al ,\be) = 0} }{
\partial \be } \right)\frac{\partial \upsilon_{(\al ,\be) = 0} }{
\partial \be } =  \frac{\,2\,}{ \omega_{i}^{*} \upsilon_i } \sum_{ k \neq i} \frac{ (\, \omega^{*}_k S_E\, \upsilon_i\,)( \omega^{*}_i S_E \,\upsilon_k\,) }{
(\lambda_i - \lambda_k)\,\omega^{*}_k \upsilon_k } \;\;\;\;\;\;~~ \;\;\\
&& \frac{\partial^{2} \lambda_{(\al ,\be) = 0} }{
\partial \al \, \partial \be }= \frac{\,1\,}{ \omega_{i}^{*} \upsilon_i }\, \left(  \omega_i^{*} H_E \frac{\partial  \upsilon_{(\al ,\be) = 0} }{
\partial \be } + \omega_i^{*} S_E  \frac{\partial  \upsilon_{(\al ,\be) = 0}  }{
\partial \al } \right)  \nonumber  \\
&& ~~~~~~ ~~~~~~~~ = \frac{\,1\,}{ \omega_{i}^{*} \upsilon_i }
\sum_{ k \neq i} \frac{ ( \omega^{*}_i H_E \,\upsilon_k\,)(\,
\omega^{*}_k S_E\, \upsilon_i\,) + (\, \omega^{*}_k H_E\,
\upsilon_i\,)( \omega^{*}_i S_E \,\upsilon_k\,) }{ (\lambda_i -
\lambda_k)\,\omega^{*}_k \upsilon_k }.~~~~~~~~~~~~~~~~~~~~~~
\nonumber
\end{eqnarray}
Therefore, the second differential $\,d^{2}\lambda_{\al ,\be}\,$ is
equal to
\newpage
\begin{eqnarray*}
d^{2}\lambda_{\al ,\be } & = & \al^{2}
\,\frac{\partial^{2}\lambda_{(\al ,\be) = 0} }{
\partial \al^{2} } +2\al \be \,\frac{\partial^{2} \lambda_{(\al ,\be) = 0}  }{
\partial \al \, \partial \be } + \be^{2} \,\frac{\partial^{2} \lambda_{(\al ,\be) = 0}  }{
\partial \be^{2} } \\
& = & \frac{\,2\al^{2} \,}{ \omega_{i}^{*} \upsilon_i } \sum_{ k
\neq i} \frac{ (\, \omega^{*}_k H_E\, \upsilon_i\,)( \omega^{*}_i
H_E \,\upsilon_k\,) }{ (\lambda_i - \lambda_k)\,\omega^{*}_k \upsilon_k } + \frac{\,2 \al
\,\be \,}{ \omega_{i}^{*} \upsilon_i } \sum_{ k \neq i} \frac{ (
\omega^{*}_i H_E \,\upsilon_k\,)(\, \omega^{*}_k S_E\,
\upsilon_i\,) + (\, \omega^{*}_k H_E\, \upsilon_i\,)( \omega^{*}_i
S_E \,\upsilon_k\,) }{ (\lambda_i -
\lambda_k )\,\omega^{*}_k \upsilon_k}\\
& & + \frac{\,2\be^{2} \,}{ \omega_{i}^{*} \upsilon_i } \sum_{ k
\neq i} \frac{ (\, \omega^{*}_k S_E\, \upsilon_i\,)( \omega^{*}_i
S_E \,\upsilon_k\,) }{ (\lambda_i - \lambda_k)\,\omega^{*}_k \upsilon_k }
 \\
& = & \frac{\,2\al \,}{ \omega_{i}^{*} \upsilon_i }\sum_{ k \neq
i} \frac{ (\, \omega^{*}_i\, H_E\, \upsilon_k\,)( \omega^{*}_k\,
\Le \,\upsilon_i\,) }{ (\lambda_i - \lambda_k)\,\omega^{*}_k \upsilon_k } + \frac{\,2\be \,}{
\omega_{i}^{*} \upsilon_i }\sum_{ k \neq i} \frac{ (\,
\omega^{*}_i\, S_E\, \upsilon_k\,)( \omega^{*}_k\, \Le
\,\upsilon_i\,) }{ (\lambda_i - \lambda_k)\,\omega^{*}_k \upsilon_k }
\\
& = & \frac{\,2 \,}{\,\omega_{i}^{*} \upsilon_i }\sum_{ k \neq i}
\frac{ (\, \omega^{*}_i\, \Le\, \upsilon_k\,)( \omega^{*}_k \,\Le
\,\upsilon_i\,) }{ (\lambda_i - \lambda_k)\,\omega^{*}_k \upsilon_k},
\end{eqnarray*}
and by
\begin{eqnarray*}
\lab = \la_i + d\lambda_{\al ,\be } + \frac{
\,1\,}{2}d^{\,2}\lambda_{\al ,\be } + \mathcal{O}\left(
\al^{3},\be^{3}\right)
\end{eqnarray*} we receive (\ref{eigen1}), whereas we have declared $\,s_{\ell}= \,\omega^{*}_{\ell} \upsilon_{\ell} .$
\newline
Multiplying the first of (\ref{deriv2}) by $\,\omega_j^{*},\,$ due
to $\,\omega_j^{*} \upsilon_i = 0\,$ $\,(j \neq i),\,$ for $\,
(\al, \be) \rightarrow ( 0, 0), \,$ we obtain
\begin{eqnarray*}\,
(\lambda_i - \lambda_j )\,\omega^{*}_j \,\frac{\partial^{2}
\upsilon_{(\al ,\be) = 0} }{
\partial \al^{2}}=  2 \left ( \omega^{*}_j H_E\,  - \omega^{*}_j \, \frac{\partial
 \lambda_{(\al ,\be)=0 }}{ \partial \al } \right)\frac{\partial \upsilon_{(\al ,\be)
= 0} }{ \partial \al }.
\end{eqnarray*}
Substituting the formulae of $\,\frac{\partial \lambda_{(\al
,\be)=0} }{
\partial \al }, \, \frac{\partial \upsilon_{(\al ,\be)=0} }{
\partial \al }\,$ from (\ref{elabda1}) and (\ref{equa1}), since $\,\omega_j^{*} \upsilon_k =
0\,$ $\,(j \neq k),\,$ we take
$$
\omega^{*}_j \,\frac{\partial^{2} \upsilon_{(\al ,\be) = 0} }{
\partial \al^{2}}=
2 \, \left( \sum_{k \neq i }\frac{ ( \omega_k^{*} H_E \,
\upsilon_i ) ( \omega_j^{*} H_E \, \upsilon_k ) }{( \lambda_i -
\lambda_k )(\lambda_i - \lambda_j)\,\omega^{*}_k \upsilon_k\,} \right) - \frac{ 2(
\omega_i^{*} H_E \, \upsilon_i ) ( \omega_j^{*} H_E \, \upsilon_i
) }{ \omega^{*}_i \upsilon_i\,( \lambda_i - \lambda_j )^{2} }
\;\;;\;\;\;\;j \neq i
$$
and then
\begin{eqnarray}
\label{equa22}
 \frac{\partial^{2} \upsilon_{(\al ,\be) = 0} }{
\partial \al^{2}} =  2
 \sum_{j \neq i} \, \left( \sum_{k \neq i }\frac{ ( \omega_k^{*} H_E \, \upsilon_i ) (
\omega_j^{*} H_E \, \upsilon_k ) }{ ( \lambda_i - \lambda_k) (
\lambda_i - \lambda_j )\,(\omega^{*}_k \upsilon_k)\,(\omega^{*}_j
\upsilon_j) }\,\right) \upsilon_j
  - 2 \sum_{j \neq i }\frac{ (
\omega_i^{*} H_E \, \upsilon_i ) ( \omega_j^{*} H_E \, \upsilon_i
) }{ \,( \lambda_i - \lambda_j )^{2}\, (\omega^{*}_i
\,\upsilon_i)(\omega^{*}_j \upsilon_j) }\,\upsilon_j .
\end{eqnarray}
Similarly, the last two expressions of (\ref{deriv2}) lead to
\begin{eqnarray}
\label{equa23}
\frac{\partial^{2} \upsilon_{(\al ,\be) = 0} }{
\partial \be^{2}}  =  2 \sum_{j \neq i }\left( \sum_{k \neq i } \frac{ ( \omega_k^{*} S_E \, \upsilon_i ) ( \omega_j^{*} S_E \, \upsilon_k )  }{ ( \lambda_i - \lambda_k)
( \lambda_i - \lambda_j )(\omega^{*}_k \upsilon_k)\,(\omega^{*}_j
\upsilon_j) }\right) \upsilon_j - 2 \sum_{j \neq i }\frac{ (
\omega_i^{*} S_E \, \upsilon_i ) ( \omega_j^{*} S_E \, \upsilon_i
) }{( \lambda_i - \lambda_j )^{2} (\omega^{*}_i \,\upsilon_i)
(\omega^{*}_j
\,\upsilon_j) }\,\upsilon_j ~~~~~\nonumber  \\
\\
\frac{\partial^{2} \upsilon_{(\al ,\be) = 0}}{
\partial \al \partial \be } = \sum_{j \neq i} \left(\sum_{k \neq i} \frac{ ( \omega_k^{*} S_E \, \upsilon_i ) ( \omega_j^{*} H_E \, \upsilon_k ) + ( \omega_k^{*} H_E \, \upsilon_i ) ( \omega_j^{*} S_E \, \upsilon_k )  }{ ( \lambda_i - \lambda_k) ( \lambda_i - \lambda_j )
(\omega^{*}_k \upsilon_k)\,(\omega^{*}_j \upsilon_j)}\right) \upsilon_j ~~~~~~~~~~~~~~~~~~~~~~~~~~~~~~~~~~~~~ \nonumber \\
~~~~~~~~~ - \sum_{j \neq i} \frac{ \left[ ( \omega_j^{*} S_E \,
\upsilon_i ) ( \omega_i^{*} H_E \, \upsilon_i ) + ( \omega_j^{*}
H_E \, \upsilon_i ) ( \omega_i^{*} S_E \, \upsilon_i )\right] }{
( \lambda_i - \lambda_j)^{2}(\omega^{*}_i\,\upsilon_i)(\omega^{*}_j \upsilon_j)
}\,\upsilon_j. ~~~~~~~~~~~~~~~~~~~~ \nonumber
\end{eqnarray}
Thus, by (\ref{equa22}) and (\ref{equa23}) we take
\begin{eqnarray*}
d^{2}\upsilon_{\al ,\be } & = & \al^{2}
\,\frac{\partial^{2}\upsilon_{(\al ,\be) = 0} }{
\partial \al^{2} } +2\al \be \,\frac{\partial^{2} \upsilon_{(\al ,\be) = 0} }{
\partial \al \, \partial \be } + \be^{2} \,\frac{\partial^{2} \upsilon_{(\al ,\be) = 0}  }{
\partial \be^{2} } \\
& = & 2 \al \sum_{j \neq i}\sum_{k \neq i} \frac{ ( \omega_j^{*}\,
H_E \, \upsilon_k ) ( \omega_k^{*}\,\Le \, \upsilon_i ) }{ (
\lambda_i - \lambda_k) ( \lambda_i - \lambda_j )(\omega^{*}_k \upsilon_k)\,(\omega^{*}_j \upsilon_j) }\,\upsilon_j + 2
\be \sum_{j \neq i}\sum_{k \neq i} \frac{ ( \omega_j^{*}\, S_E \,
\upsilon_k ) ( \omega_k^{*}\,\Le \, \upsilon_i ) }{ ( \lambda_i -
\lambda_k) (
\lambda_i - \lambda_j )(\omega^{*}_k \upsilon_k)\,(\omega^{*}_j \upsilon_j) }\,\upsilon_j \\
& & - \frac{\,2 \al \,}{\omega^{*}_i\,\upsilon_i }\sum_{j \neq i}
\frac{ ( \omega_i^{*}\, H_E \, \upsilon_i ) ( \omega_j^{*}\,\Le \,
\upsilon_i )}{ ( \lambda_i - \lambda_j)^{2}(\omega^{*}_j \upsilon_j)  }\,\upsilon_j -
\frac{\,2 \be \,}{\omega^{*}_i\,\upsilon_i }\sum_{j \neq i} \frac{
( \omega_i^{*}\, S_E \, \upsilon_i ) ( \omega_j^{*}\,\Le \,
\upsilon_i ) }{ ( \lambda_i - \lambda_j)^{2}(\omega^{*}_j \upsilon_j)  }\,\upsilon_j\\
& = & 2 \sum_{j \neq i} \sum_{k \neq i} \frac{ ( \omega_j^{*}\,
\Le \, \upsilon_k ) ( \omega_k^{*}\,\Le \, \upsilon_i ) }{ (
\lambda_i - \lambda_k) ( \lambda_i - \lambda_j )(\omega^{*}_k \upsilon_k)\,(\omega^{*}_j \upsilon_j) }\,\upsilon_j \\
&& ~~~~~~~~~~~~~~~~~~~~~~~~~~~~~~~~~ - \frac{\,2
\,}{\omega^{*}_i\,\upsilon_i }\sum_{j \neq i} \frac{ (
\omega_i^{*}\, \Le \, \upsilon_i ) ( \omega_j^{*}\,\Le \,
\upsilon_i )}{ ( \lambda_i - \lambda_j)^{2} (\omega^{*}_j \upsilon_j) }\,\upsilon_j ,
\end{eqnarray*}
and by
$$  \uab=  \upsilon_i +d\upsilon_{\al ,\be}+\frac{\,1\,}{2} d^{\,2}\upsilon_{\al ,\be } +\mathcal{O}
\left( \al^{3},\be^{3}\right) $$  we obtain the claimed equality
(\ref{eigen2}).
 \qed  \ \\

The simplified presentation of partial differential formulae
(\ref{elabda1}) and (\ref{eigall}) of $\lab\,$ and (\ref{equa1}),
(\ref{equa2}), (\ref{equa22}) and (\ref{equa23}) of $\uab\,$ for
$\,\al=\beta =0 \,$ are independent results on those, which were
obtained earlier by Chu in \cite{Chu}.
 Chu has follow different methodology considering that
$\,\la_i\,$ is simple eigenvalue and an additional normalized
condition that $\,\omega_i^{*} \upsilon_i =1,\,$ and even Chu's
formulations of the partial derivatives depend on the invertibility
of a matrix and the eigenvectors $\,\upsilon_i,\,\omega_i.\,$
Furthermore, no results on the perturbation of the eigenpairs
$\,\lab, $ and $\,\uab\,$ are given in \cite{Chu}.

 In the following we present a
lemma, which will contribute in the approximation formulae
(\ref{eigen1}) and (\ref{eigen2}).
\vspace{0.1cm}

\begin{lema}
\label{lem1} ~Let the matrix $\,A \in \cM\,$ be diagonalizable and
$\,Y_i,\,$ $\,W_i\,$ be matrices whose columns $\,\upsilon_i\,$
and rows $\,\omega_i^{*}\,$ respectively are the corresponding
right and left eigenvectors of $\,A,\,$ for $\,\lambda_i \in
\sigma(A).\,$ A generalized inverse of $\,(A-\lambda_i I)^{\mu},
\,$ $\,\mu \in \mathbb{N},\,$ is defined by
\begin{eqnarray}
\label{genarinv} \left[(A-\lambda_i I )^{\mu}\right]^{+}\,=\,
\sum_{k \neq i\,} \frac{ \upsilon_k \,\omega_k^{*}\,}{ (\lambda_k
- \lambda_i)^{\mu}\,s_k} \;\;;\;\;s_k=\omega^{*}_k \upsilon_k .
\end{eqnarray}
\end{lema}
\pr~It is evident that $\,(A-\lambda_i I )\frac{ \upsilon_k
\,\omega_k^{*}\,}{ \lambda_k - \lambda_i} = \upsilon_k
\,\omega_k^{*},\,$ and then
\begin{eqnarray*}
\left( A-\lambda_i I \right)^{\mu} \frac{ \upsilon_k
\,\omega_k^{*}\,}{( \lambda_k - \lambda_i)^{\mu}} & = & \left(
A-\lambda_i I \right)^{\mu-1}\left( A-\lambda_i I \right)\frac{
\upsilon_k \,\omega_k^{*}\,}{( \lambda_k - \lambda_i)^{\mu}}  \\
& = &  \left( A-\lambda_i I \right)^{\mu-1} \frac{ \upsilon_k
\,\omega_k^{*}\,}{( \lambda_k - \lambda_i)^{\mu-1}} = \cdots =
\upsilon_k\,\omega^{*}_k.
\end{eqnarray*}
Since $\,\sum \frac{\upsilon_k \,\omega_k^{*}}{\,s_k\,} = I,\,$
$\,\omega_i^{*}A = \lambda_i \omega_i^{*}\,$ and for $\, k \neq
i,\,$ $\,\omega_k^{*} \,\upsilon_i=0, \,$ we have:
\begin{eqnarray*}
\left( A-\lambda_i I \right)^{\mu} \left( \sum_{k \neq i\,} \frac{
\upsilon_k \,\omega_k^{*}\,}{ (\lambda_k -
\lambda_i)^{\mu}\,s_k}\right)\left( A-\lambda_i I \right)^{\mu} =
\left( \sum_{k \neq i\,} \frac{\upsilon_k \,\omega_k^{*}}{\,s_k\,} \right)\left(
A-\lambda_i I \right)^{\mu}=\left( I- Y_i W_i \right)
\left( A-\lambda_i I
\right)^{\mu} ~~~~~~ \\
= (A-\lambda_i I)^{\mu}- Y_iW_i(A-\lambda_i I)^{\mu} =
(A-\lambda_i I)^{\mu}-Y_i \left(W_iA - \lambda_i W_i
\right)(A-\lambda_i I)^{\mu-1}=(A-\lambda_i I)^{\mu}
\end{eqnarray*}
and
\begin{eqnarray*}
\left( \sum_{k \neq i\,} \frac{ \upsilon_k \,\omega_k^{*}\,}{
(\lambda_k - \lambda_i)^{\mu}\,s_k }\right) \left( A-\lambda_i I
\right)^{\mu} \left( \sum_{k \neq i\,} \frac{ \upsilon_k
\,\omega_k^{*}\,}{ (\lambda_k - \lambda_i)^{\mu}\,s_k }\right) ~~~~~~~~~~~~~~~~~~~~~~~~~~~~~~\\
~~~~~~~~~~~~~~~~~= \left(
\sum_{k \neq i\,} \frac{ \upsilon_k \,\omega_k^{*}\,}{ (\lambda_k
- \lambda_i)^{\mu}\,s_k}\right) \left( I - Y_i \,W_i \right) =
\sum_{k \neq i\,} \frac{ \upsilon_k \,\omega_k^{*}\,}{ (\lambda_k
- \lambda_i)^{\mu}\,s_k}.
\end{eqnarray*}
\qed

 In Lemma \ref{lem1}, if $\,A\,$ is normal, then $\,\upsilon_k =
\omega_k,\,$ and $\,\left[(A-\lambda_i I )^{\mu}\right]^{+}\,$ is
Hermitian. In this case, we confirm that $\, \left[(A-\lambda_i I
)^{\mu}\right]^{+}\,$ in (\ref{genarinv}) is the {\it Moore-Penrose
inverse} of $\, (A-\lambda_i I )^{\mu}.\,$ \vspace{1.2cm}

Combining Equations (\ref{genarinv}),  (\ref{eigen1}), and
(\ref{eigen2}), Theorem \ref{normaleig} leads to a generalization of
a corresponding result for simple eigenvalues of a Hermitian matrix,
which was presented by M. Fiedler in \cite{Fi}.

 \vspace{0.1cm}

\begin{them}
\label{eiggeneral} ~Let the matrix $\,A \in \cM\,$ be diagonalizable
and $\, \lambda_i\,$ be a semisimple eigenvalue of $\,A\,$ with
$\,\upsilon_i,\,\omega_i\,$ corresponding right and left
eigenvectors. If the assumptions for the equations (\ref{eigen1})
and (\ref{eigen2}) hold, then the following expressions for
$\,\lab\,$ and $\,\uab\,$ hold:
\begin{eqnarray}
\label{egen1} \lab & = & \la_i + \frac{\,1\,}{\,s_i\,}\omega^{*}_i \Le \upsilon_i
-\frac{\,1\,}{\,s_i\,}\omega_{i}^{*}\Le  (A-\lambda_i I )^{+} \Le \,\upsilon_i +
\mathcal{O}\left( \al^{3},\be^{3}\right),~~~~~ \end{eqnarray}
\begin{eqnarray}\label{egen2} \uab & = & \upsilon_i - (A-\lambda_i I
)^{+}\, \Le \,\upsilon_i +
\left[(A-\lambda_i I )^{+}\,\Le \right]^{2}\,\upsilon_i  \nonumber \\
& & ~~~~~~~~~~~~~~~~~-\frac{\,1\,}{\,s_i\,}\, \left[(A-\lambda_i I )^{2}\right]^{+}\,\Le
\,\upsilon_i \omega_i^{*}\,\Le \,\upsilon_i + \mathcal{O}\left(
\al^{3},\be^{3}\right).~~~~~
\end{eqnarray}
\end{them}
\pr ~ From (\ref{eigen1}) and (\ref{genarinv}) with $\,\mu =1,\,$ we
immediately  have
\begin{eqnarray*}
\label{eigengen3} \lab = \la_i + \frac{\,1\,}{\,s_i\,}\,\omega^{*}_i \Le \upsilon_i
-\frac{\,1\,}{\,s_i\,}\,\omega_{i}^{*}\Le \, \left(\sum_{k \neq i\,} \frac{ \upsilon_k
\,\omega_k^{*}\,}{\, (\lambda_k - \lambda_i)\,s_k} \right) \, \Le
\,\upsilon_i\,+ \mathcal{O}\left( \al^{3},\be^{3}\right),
\end{eqnarray*}
proving (\ref{egen1}). Also, from (\ref{eigen2}) and
(\ref{genarinv}) with $\,\mu=1,2,\,$ we have
\begin{eqnarray*}
\label{eigengen4} \uab & =& \upsilon_i -\left( \sum_{k \neq i}
\frac{ \,\upsilon_k \, \omega^{*}_k }{(\lambda_k - \lambda_i)\,s_k
 }\right) \Le \,\upsilon_i
 + \sum_{k \neq i} \left(  \sum_{j \neq i}\frac{  \upsilon_j\,\omega_j^{*} }{( \lambda_j - \lambda_i)\,s_j
 }\Le \,
 \frac{  \upsilon_k \omega_k^{*} }{  (\lambda_k - \lambda_i)\,s_k }\,\Le \,\upsilon_i \right)
 \\
& & ~~~~~~~~~~~~~~~~~~ - \left( \sum_{j \neq i }\frac{ \upsilon_j \,\omega_j^{*}
}{ ( \lambda_j - \lambda_i )^{2}s_j }\, \,\right) \Le \,\frac{\upsilon_i\,
\omega_i^{*}}{\,s_i\,} \,\Le\,\upsilon_i + \mathcal{O}\left(
\al^{3},\be^{3}\right) \\
&= & \upsilon_i - (A-\lambda_i I )^{+}\, \Le \,\upsilon_i +
(A-\lambda_i I )^{+}\,\Le\,\left( \sum_{k \neq i }\frac{
\upsilon_k \,\omega_k^{*} }{ (\lambda_k - \lambda_i)\,s_k } \Le \,
\upsilon_i \,\right) \\
& &  ~~~~~~~~~~~~~~~~~~- \left( \sum_{j \neq i }\frac{ \upsilon_j \,\omega_j^{*} }{ (
\lambda_j - \lambda_i )^{2} s_j}\, \,\right) \Le \,\frac{\upsilon_i\,
\omega_i^{*}}{\,s_i\,} \,\Le\,\upsilon_i + \mathcal{O}\left(
\al^{3},\be^{3}\right)
\\
&= & \upsilon_i - (A-\lambda_i I )^{+}\, \Le \,\upsilon_i +
\left[(A-\lambda_i I )^{+}\,\Le \right]^{2}\,\upsilon_i \\
&& ~~~~~~~~~~~~~~~~~~~~~~~~~ - \left[(A-\lambda_i I
)^{2}\right]^{+}\,\Le \, \frac{\upsilon_i
\omega^{*}_i}{\,s_i\,}\,\Le \,\upsilon_i + \mathcal{O}\left(
\al^{3},\be^{3}\right),
\end{eqnarray*}
proving (\ref{egen2}).
 \qed \vspace{0.2cm}
\newline
In (\ref{egen1}) and (\ref{egen2}), if we consider the first-order
approximation, then we simply have
\begin{eqnarray}
\label{e26}
\widetilde{\lambda}_{\al,\be}= \lambda_i+
\frac{\,1\,}{\,s_i\,}\omega^{*}_i
\Le \, \upsilon_i + \mathcal{O}\left( \al^{2},\be^{2}\right),~~\hspace{2.9cm}~~~~~~~ \nonumber \\
\mbox{and} ~~~~~~~~~~~~~~~~~~~~~~~~~~~~~~~~~~~\hspace{6cm} ~~~~~~~~~~~~~~~~~~~~~~~~~~~\\
 \widetilde{\upsilon}_{\al,\be}=
\upsilon_i-(A- \la_i I )^{+} \Le\,\upsilon_i + \mathcal{O}\left(
\al^{2},\be^{2}\right), ~~\hspace{2cm}~~~~~~~\nonumber
\end{eqnarray}
for a simple eigenvalue $\,\lambda_i.\,$ In these cases,
\begin{eqnarray*}
M_{\al,\be}\widetilde{\upsilon}_{\al,\be}-\widetilde{\la}_{\al,\be}\widetilde{\upsilon}_{\al,\be}& = & -(A-\la_i
I)(A-\la_i I)^{+} \Le \, \upsilon_i + \Le \, \upsilon_i - \,\frac{\,1\,}{s_i\,}\left[
\omega_i^{*} \, \Le \, \upsilon_i \right] \upsilon_i \\
& & ~~~~~ +\left[\,\frac{\,1\,}{s_i\,}(\omega_i^{*} \, \Le \,\upsilon_i)I  - \Le \right](A-\la_i I)^{+} \Le \, \upsilon_i   \\
& = & -\left[(A-\la_i I)(A-\la_i I)^{+} + \frac{\upsilon_i \omega_i^{*}}{\,s_i\,}
\right]
\, \Le\, \upsilon_i + \Le\, \upsilon_i \\
& & ~~~~ +\left[ \,\frac{\,1\,}{s_i\,}(\omega_i^{*} \, \Le\,\upsilon_i ) I - \Le \right](A-\la_i
I)^{+} \Le \, \upsilon_i +  \mathcal{O}\left( \al^{2},\be^{2}\right).
\end{eqnarray*}
Since,
$$
(A-\la_i I)(A-\la_i I)^{+}\,+\, \frac{\upsilon_i \omega_i^{*}}{\,s_i\,}= (A-\la_i I)\sum_{ k \neq i}\frac{
\upsilon_k \,\omega_k^{*} }{ \,(\lambda_k - \lambda_i)\,s_k }\,+ \,\frac{\upsilon_i \omega_i^{*}}{\,s_i} = \sum_{k} \frac{\upsilon_k \omega_k^{*}}{\,s_k\,}=I,
$$
we have
\begin{eqnarray}
\label{egen3}
M_{\al,\be}\widetilde{\upsilon}_{\al,\be}-\widetilde{\la}_{\al,\be}\widetilde{\upsilon}_{\al,\be}
=\left[ \,\frac{\,1\,}{\,s_i\,}\omega_i^{*} \, \Le \, \upsilon_i I -
\Le \right] (A-\la_i I)^{+} \Le  \upsilon_i +\mathcal{O}\left(
\al^{2},\be^{2}\right).
\end{eqnarray}

\begin{prop}
\label{eiggll} ~Let $\,\lambda_i \,$ be a simple eigenvalue of
diagonalizable matrix $\,\,A \in \cM \,\,$ with right and left
eigenvectors $\,\upsilon_i \,$ and $\,\omega_i.\,$ If there exist
$\,\al,\be\,$ such that $\,\,\Le\,\upsilon_i \in ker(A-\la_i
I)^{+},\,$ then $\,\widetilde{\la}_{\al,\be}\,$ and
$\,\widetilde{\upsilon}_{\al,\be}\,$ in (\ref{e26}) is an
approximation of an eigenpair of $\,M_{\al,\be}=A + \Le.$
\end{prop}
\begin{coro}
\label{eigcor} ~Let $\,\lambda_i \,$ be a simple eigenvalue of
normal matrix $\,A \in \cM\,$ with eigenvector $\,\upsilon_i. \,$ If
there exist $\, \al, \be \,$ such that $\,\Le\,\upsilon_i
=\overline{\upsilon}_i ,\, $ then $\,\widetilde{\la}_{\al,\be}\,$
and $\,\widetilde{\upsilon}_{\al,\be}\,$ in (\ref{e26}) is an
approximation of an eigenpair of $\,M_{\al,\be}=A + \Le.$
\end{coro}
\pr~ It is well-known that
$$
ker(A-\la_i I)^{+}= ker(A-\la_i I)^{T}.
$$
Also, since $\,A\,$ is normal, $\,\omega_i = \upsilon_i.\,$ Thus by
(\ref{egen3}), it is implied that
$$
(A-\la_i I)^{+} \Le \, \upsilon_i = (A-\la_i I)^{+}\,\overline{\upsilon_i} = \left[ \upsilon_i^{*}(A-\la_i I) \right]^{T} = 0.
$$
\qed

In this case, (\ref{e26}) is simplified to
\begin{eqnarray}
\label{eghen4} \widetilde{\lambda}_{\al,\be}= \lambda_i +
\frac{\,\upsilon_i^{*} \overline{\upsilon_i}\,}{\,s_i\,}+
\mathcal{O}\left( \al^{2},\be^{2}\right),\;\;\;\;\;\;\;\;
\widetilde{\upsilon}_{\al,\be}= \upsilon_i+\mathcal{O}\left(
\al^{2},\be^{2}\right)
\end{eqnarray}
or, $\,\,\widetilde{\lambda}_{\al,\be}= (\lambda_i +1)+
\mathcal{O}\left( \al^{2},\be^{2}\right),\;\;\;
\widetilde{\upsilon}_{\al,\be}= \upsilon_i+\mathcal{O}\left(
\al^{2},\be^{2}\right)\,\,\,$ for real symmetric matrix $\,A.\,$
\vspace{0.8cm}
\newline
{\bf Example 1}~ Let $\,A=\left[ \begin{array}{ccc}
  1 & 0 & 1 \\
  0 & 2 & 0 \\
  0 & 0 & 2 \\
\end{array}\right].\,
\,$ Then $\,\, \sigma(A) =\{ \,\la_1=1,\,\,\la_{2,3}=2 \, \},\,\,$
and the corresponding right and left eigenvectors are given by :
$$\,\, \,\left[
\begin{array}{ccc}
  \upsilon_1 & \upsilon_2 & \upsilon_3
\end{array}
\right] = \left[ \begin{array}{ccc}
  1 & 1 & 0 \\
  0 & 0 & 1 \\
  0 & 1 & 0 \\
\end{array}\right], \,\;\;\;\;\;\,\left[
\begin{array}{ccc}
  \omega_1 & \omega_2 & \omega_3
\end{array}
\right] = \left[ \begin{array}{ccc}
  1 & 0 & 0 \\
  0 & 0 & 1 \\
  -1 & 1 & 0 \\
\end{array}\right].\,\,$$
\newline
Let the perturbation matrix $\,E=\left[
\begin{array}{ccc}
  1 & 0 & 1 \\
  0 & 1 & 0 \\
  -1 & 0 & 0 \\
\end{array}\right].
\,$ We then obtain $$\,{\cal L}(E,\al, \beta)= \left[
\begin{array}{ccc}
  \al & 0 & \be \\
  0 & \al & 0 \\
  -\be & 0 & 0 \\
\end{array}\right]\,\,\;\; \mbox{ and }\,\;\;\,M_{\al, \be}=\left[ \begin{array}{ccc}
  1+\al & 0 & 1+\be \\
  0 & 2+\al & 0 \\
  -\be & 0 & 2 \\
\end{array}\right].
\,$$ Then,  $\, \sigma(M_{\al, \be}) =\{ \lambda_{\small
{1,\al,\be}}=\frac{3+\al -\sqrt{(\al-1)^{2}-4\beta
(\beta+1)}\,}{\,2\,},\;\,\lambda_{\small{2,\al,\be}}=\frac{3+\al
+\sqrt{(\al-1)^{2}-4\beta (\beta+1)}\,}{\,2\,}, \,\;\,\newline
\;\;\lambda_{\small{3,\al,\be}}=2+\al  \},\,$ and the corresponding
right eigenvectors are:
$$\,\left[\begin{array}{ccc} 1+\be & 0 & \lambda_{1,\al,\be}-1-\al
\end{array}\right]^{T},\,\;\;\;\,\,\left[\begin{array}{ccc} 1+\beta & 0 &
\lambda_{2,\al,\be}-1-\al
\end{array}\right]^{T},\;\;\;\,\,\,\left[\begin{array}{ccc} 0 &
1 & 0 \end{array}\right]^{T}.$$ Clearly, the eigenvalues and the
corresponding eigenvectors of $\,M_{\al, \be}\,$ are real functions
of two variables, with continuous partial derivatives for all
permissible values of $\,\al,\be.$
\newline For $\,\al=0.1,\, \be=0.01,\,$ we have $\, \sigma(M_{\al,
\be}) =\{ \,1.1114,\,\,1.9886,\,\,2.1 \, \},\,$ and the
corresponding unit eigenvectors are given by
$$\,\left[\begin{array}{ccc} 0.9999 & 0 & 0.0113
\end{array}\right]^{T},\,\;\;\;\,\,\left[\begin{array}{ccc} 0.7508 & 0 &
0.6606
\end{array}\right]^{T},\;\;\;\,\,\,\left[\begin{array}{ccc} 0 &
1 & 0 \end{array}\right]^{T}.$$ Moreover by (\ref{genarinv}),
$$\, (A-I)^{+} = [(A-I)^{2}]^{+}=\left[
\begin{array}{ccc}
  0 & 0 & 1 \\
  0 & 1 & 0 \\
  0 & 0 & 1 \\
\end{array}\right]. \,$$
Hence, a first-order approximations of the eigenvalue and unit
eigenvector of $\,M_{\al,\be}\,$ by (\ref{e26}) are:
$$
\widetilde{\lambda}_1=1+0.11=1.11,
\;\;\;\;\widetilde{\upsilon}_1^{T}=\left[\begin{array}{ccc} 1 & 0 &
0.0099 \end{array}\right]^{T}.
$$
Also, by (\ref{egen1}), (\ref{egen2}), the corresponding
second-order approximations are equal to
$$
\widetilde{\lambda}_1=1.1112,\;\;\;\;\;\;
\widetilde{\upsilon}_1^{T}=\left[\begin{array}{ccc} 0.9999 & 0 &
0.0111 \end{array}\right]^{T}.
$$
By (\ref{genarinv}), $\, (A-2I)^{+} = \left[
\begin{array}{ccc}
  -1 & 0 & 1 \\
  0 & 0 & 0 \\
  0 & 0 & 0 \\
\end{array}\right] , \,$ $\, [(A-2I)^{2}]^{+} = \left[ \begin{array}{ccc}
  1 & 0 & -1 \\
  0 & 0 & 0 \\
  0 & 0 & 0 \\
\end{array}\right] , \,$
and substituting these expressions in (\ref{egen1}) and
(\ref{egen2}), using the eigenvectors  $\,\upsilon_2, \omega_2,\,$
the second-order approximations of the eigenpair of $\,M_{\al,
\be}\,$ are
$$
\widetilde{\lambda}_2= 1.9888,\;\;\;\;
\widetilde{\upsilon}_2^{T}=\left[\begin{array}{ccc} 0.7505 & 0 &
0.6609 \end{array}\right]^{T}.
$$
Similarly, using the eigenvectors $\,\upsilon_3, \omega_3\,$, the
second-order approximations of the third eigenpair of $\,M_{\al,
\be}\,$ are
$$
\widetilde{\lambda}_3 = 2.1,\;\;\;\;
\widetilde{\upsilon}_3^{T}=\left[\begin{array}{ccc} 0 & 1 & 0
\end{array}\right]^{T}.
$$
Notice that, we obtain all eigenpair of $\,M_{\al, \be}\,$ with
preciseness $10^{-3}$. \newline {\bf Example 2}~~Let the matrix
$\,A=\left[ \begin{array}{ccc}
  2 & -2 & 0 \\
  -2 & 6 & -1 \\
  0 & -1 & 2 \\
\end{array}\right],
\,$ with $\, \sigma(A) =\{ \,\la_1=2,\,\,\la_{2}=1,\,\,\la_3=7 \,
\},\,$  and $\,\left[
\begin{array}{ccc}
  \upsilon_1 & \upsilon_2 & \upsilon_3
\end{array}
\right] = \left[ \begin{array}{ccc}
  1 & 2 & 2 \\
  0 & 1 & -5 \\
  -2 & 1 & 1 \\
\end{array}\right]\,$ be the corresponding right and left eigenvectors.
\newline
For $\,\al=0.04,\, \be=0.08,\,$ and non symmetric $\,E=\left[
\begin{array}{ccc}
  21 & 3 & -4 \\
  1 & 1 & 0 \\
  -8 & 0 & 20 \\
\end{array}\right],
\,$ then
\newline
$\, {\cal L}(E,0.04, 0.08)= \left[ \begin{array}{ccc}
  0.84 & 0.16 & -0.08 \\
  0 & 0.04 & 0 \\
  -0.4 & 0 & 0.8 \\
\end{array}\right]\,\,\, $  and
$\,\,\,M_{\al, \be}=\left[ \begin{array}{ccc}
  2.84 & -1.84 & -0.08 \\
 -2 & 6.04 & -1 \\
  -0.4 & -1 & 2.8 \\
\end{array}\right].
\,$ \newline We have  $\, \sigma(M_{\al, \be}) =\{
\,3,\,\,1.5886,\,\,7.0914 \, \},\,$ with corresponding eigenvectors
$$\,\left[\begin{array}{ccc}
1 & 0 & -2
\end{array}\right]^{T},\,\;\,\left[\begin{array}{ccc}
0.6814 & 0.4379 & 0.5865
\end{array}\right]^{T},\,\;\,\left[\begin{array}{ccc} 0.3883 &
-0.9048 & 0.1747 \end{array}\right]^{T}.\,$$ Since $\, {\cal
L}(E,0.04, 0.08)\upsilon_1 = \upsilon_1,\,$ by Corollary
\ref{eigcor} and (\ref{eghen4}) the eigenpair of $\,M_{\al,\be}\,$
is
$$
\widetilde{\lambda}_1=\lambda_1+\frac{\upsilon_1^{*}
\upsilon_1}{s_1}= 2 + 1= 3,
\;\;\;\;\;\widetilde{\upsilon}_1^{T}=\upsilon_1^{T}.
$$
Moreover, $\, {\cal L}(E,0.04, 0.08)\upsilon_2 \neq  \upsilon_2,\,$
and by (\ref{genarinv}), $\, (A-I)^{+} = \left[ \begin{array}{ccc}
  0.222 & -0.0556 & -0.3889\\
  -0.0556 & 0.1389 & -0.0278 \\
  -0.3889 & -0.0278 & 0.8056 \\
\end{array}\right] , \,$ $\, [(A-I)^{2}]^{+} = \left[ \begin{array}{ccc}
  0.2037 & -0.0093 & -0.3981 \\
  -0.0093 & 0.0231 & -0.0046 \\
  -0.3981 & -0.0046 & 0.8009 \\
\end{array}\right] .\,$  Then, (\ref{egen1}) and (\ref{egen2}) lead to the second-order
approximations
$$
\widetilde{\lambda}_2=1.5890,\;\;\;\;\;
\widetilde{\upsilon}_2^{T}=\left[\begin{array}{ccc} 0.7031 & 0.4409
& 0.5579 \end{array}\right]^{T},
$$
in contrast to the first-order approximations of the eigenpair of
$\,M_{\al,\be},\,$ which are equal to
$$
\widetilde{\lambda}_2=1.5933,\;\;
\widetilde{\upsilon}_2^{T}=\left[\begin{array}{ccc} 0.6257 & 0.4242
& 0.6546 \end{array}\right]^{T}.
$$
Moreover, $\, {\cal L}(E,0.04, 0.08)\upsilon_3 \neq  \upsilon_3,\,$
and by (\ref{genarinv}), $\, (A-7I)^{+} = \left[ \begin{array}{ccc}
  -0.1511 & 0.0556 & 0.0244 \\
  -0.0556 & -0.0278 & -0.0278 \\
  0.0244 & -0.0278 & -0.1878 \\
\end{array}\right] , \,$ $\, [(A-7I)^{2}]^{+} = \left[ \begin{array}{ccc}
  0.0265 & 0.0093 & -0.0067 \\
  0.0093 & 0.0046 & 0.0046 \\
  -0.0067 & 0.0046 & 0.0366 \\
\end{array}\right] . \,$  Also, by (\ref{egen1}) and (\ref{egen2}), the second-order approximations give
$$
\widetilde{\lambda}_3=7.0910,\;\;\;\;
\widetilde{\upsilon}_3^{T}=\left[\begin{array}{ccc} 0.3878 & -0.9049
& 0.1754 \end{array}\right]^{T},
$$
but the first-order approximations are not satisfactory enough
$$
\widetilde{\lambda}_3=7.0867,\;\;\;\;
\widetilde{\upsilon}_3^{T}=\left[\begin{array}{ccc} 0.3851 & -0.9056
& 0.1779 \end{array}\right]^{T}.
$$

\end{document}